\documentclass[a4paper]{amsart}

\usepackage{amssymb,amsmath,epsfig,graphics,psfrag}
\usepackage[english]{babel}

\newtheorem{theorem}{Theorem}[section]
\newtheorem{proposition}{Proposition}[section]
\newtheorem{lemma}[theorem]{Lemma}

\theoremstyle{definition}
\newtheorem{definition}[theorem]{Definition}

\theoremstyle{remark}
\newtheorem{remark}[theorem]{Remark}

\numberwithin{equation}{section}

  \def\CC{{\mathbf{C}}}

  \def\RR{{\mathbf{R}}}

 \def\ZZ{{\mathbf{Z}}}

  \def\cC{{\mathcal{C}}}
\def\cD{{\mathcal{D}}}  \def\cF{{\mathcal{F}}}
\def\cG{{\mathcal{G}}}  
  \def\cL{{\mathcal{L}}}
  
  \def\cR{{\mathcal{R}}}
 \def\cT{{\mathcal{T}}} \def\cU{{\mathcal{U}}}
\def\cV{{\mathcal{V}}}

\def\itt{{\it {(TT)}}}

\begin{document}
\sloppy

\title[Saddle connections in dimension 3]{Classification of a generic set of saddle connections in dimension 3}

\author{Emmanuel Dufraine}
\address{Mathematics Institute, University of Warwick, CV4 7AL Coventry, U.K.}
\email{dufraine@maths.warwick.ac.uk}



\begin{abstract}
We consider an open and dense subset of the set of all connections (heteroclinic orbits) between two saddles of a vector field in dimension three. We give a complete classification up to topological equivalence of the dynamics in a neighbourhood of a connection in this subset.
\end{abstract}

\maketitle


\section{Introduction, presentation of the results}

We are interested in {\em topological equivalence} of flows of smooth vector fields. We recall that two flows are topologically equivalent if there exists a homeomorphism that maps orbits to orbits, preserving the orientation given by the associated vector fields.

The classification of flows up to topological equivalence is usually possible if we have simple dynamics (without reccurence).

For example, the work of Fleitas in~\cite{fl}, gives a complete classification of gradient-like vector fields (Morse-Smale without periodic orbits) on compact three-dimensional manifolds. Recent work by Prishlyak, \cite{pr}, gives a generalisation by adding some periodic orbits. Similarly, a series of papers~\cite{bogr,bogrmepe2,bogrpe1,bola}, culminating in~\cite{bogrmepe3}, gives a classification of the more complicated situation of gradient-like diffeomorphisms on 3-manifolds. 

A {\em connection} (or {\em heteroclinic orbit}) between two saddles $p$ and $q$ is the adherence of an orbit accumulating on $p$ in negative time and on $q$ in positive time. The only connections existing in gradient-like vector fields are structurally stable. Another natural way to generalise the work of Fleitas, would be to include non structurally stable vector fields, by permitting more general saddle connections.

Before the global classification, one should address the problem of semi-local classification: classification, up to topological equivalence, in a neighbourhood of a connection. This is the purpose of this paper. Loosely, we will speak of ``the classification of connections" up to topological equivalence. 

We denote by $\cC$ the set of vector fields defined in a neighbourhood of a saddle connection $\Gamma$ between $p$ and $q$ in dimension three. By convention, $p$ and $q$ are named such that $\Gamma=W^u(p)\cap W^s(q)$. In this paper, we give a classification of connections in an open and dense subset of $\cC$.

\begin{figure}
   \begin{minipage}[b]{.49\linewidth}
    \psfrag{wsp}{$W^s_{loc}(p)$}\psfrag{wuq}{$W^u_{loc}(q)$}
    \psfrag{g}{$\Gamma$}
    \centerline{\includegraphics{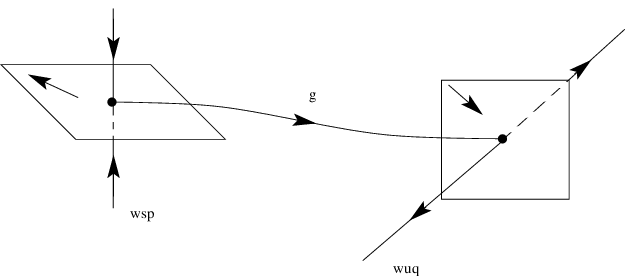}}
    \caption{$\Gamma \in$(1-2)\label{f.1-2}}  
   \end{minipage} \hfill
   \begin{minipage}[b]{.49\linewidth}
   \psfrag{wsp}{$W^s_{loc}(p)$}\psfrag{wuq}{$W^u_{loc}(q)$}
    \psfrag{g}{$\Gamma$}
    \centerline{\includegraphics{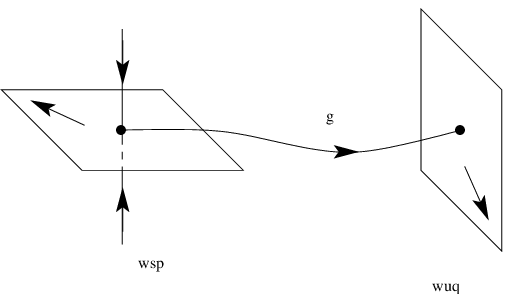}}
    \caption{$\Gamma \in$(1-1)\label{f.1-1}}
   \end{minipage}
\end{figure}

\begin{figure}
   \begin{minipage}[b]{.49\linewidth}
    \psfrag{wsp}{$W^s_{loc}(p)$}\psfrag{wuq}{$W^u_{loc}(q)$}
    \psfrag{g}{$\Gamma$}
    \centerline{\includegraphics{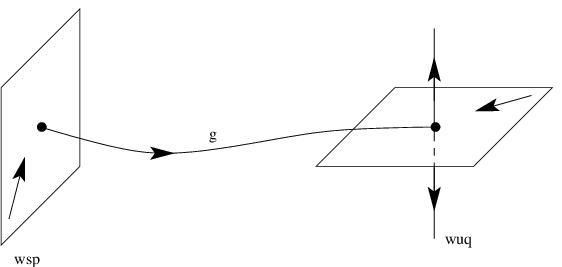}}
    \caption{$\Gamma \in$(2-2)\label{f.2-2}}  
   \end{minipage} \hfill
   \begin{minipage}[b]{.49\linewidth}
   \psfrag{wsp}{$W^s_{loc}(p)$}\psfrag{wuq}{$W^u_{loc}(q)$}
    \psfrag{g}{$\Gamma$}
    \centerline{\includegraphics{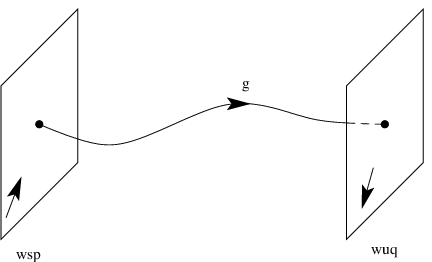}}
    \caption{$\Gamma \in$(2-1)\label{f.2-1}}
   \end{minipage}
\end{figure}

Locally, in a neighbourhood of a saddle, the Theorem of Hartman-Grobman asserts that the only topological invariant is the {\em Morse index} of the saddle (number of stable eigenvalues). This divides $\cC$ into the following subsets: $(1-2)$, $(1-1)$, $(2-2)$ and $(2-1)$, according to the Morse index of the saddles $p$ and $q$ respectively (in other words, subsets are denoted by $(\dim W^s(p),\dim W^s(q))$, compare Figures~(\ref{f.1-2}-\ref{f.2-1})). If two connections are in different subsets, they are not topologically equivalent.


In~\cite{pa}, Palis proves that in dimension two, any two saddle connections are topologically equivalent.

Many authors worked on the existence of {\em modulus of stability} for connections: topological invariants for nearby vector fields. In the paper~\cite{nepata}, Newhouse, Palis \& Takens prove the existence of a modulus of stability for a bifurcation point of a generic one-parameter family of diffeomorphisms. Takens gives a complete classification for codimension one saddle connections of vector fields in dimension 4 in~\cite{ta1,ta2}. Beloqui, \cite{belo}, started the study of vector fields in dimension three considering a connection between a saddle and a periodic orbit. 

In~\cite{vs}, van Strien generalises the works of Newhouse, Palis, Takens and Beloqui and gives a complete classification of codimension one saddle connections in any dimension. In his work, only nearby connections were considered. In particular, he proves that two nearby connections in $(1-1)$ (resp. $(2-2)$) are topologically equivalent. 

Some codimension two phenomena are studied by Vegter in~\cite{ve1,ve2} and by Dias-Carneiro \& Palis in~\cite{dipa} for families of gradient vector fields. Among other things, they prove that there is no modulus of stability for connections in $(2-1)$ with real eigenvalues. 

We prove here that in each subset $(1-1)$, $(2-2)$ and $(2-1)$ with real eigenvalues, there are two classes of topological equivalence. Let us mention that on the contrary of those previous works, we don't pay attention to the bifurcations of the connections in this paper. The concept of ``codimension" is not important in the rest of the paper: we study a generic subset of all saddle connections and not generic one or two parameters families of vector fields.

With Bonatti in~\cite{bodu}, we give a complete classification of connections in $(2-1)$ with complex eigenvalues (stable eigenvalues in $p$ and unstable eigenvalues in $q$). The surprising result in this situation is that the transition map between transverse discs in neighbourhoods of the saddles plays a key role. Indeed, according to an explicit condition between the eigenvalues at $p$ and $q$ and the modulus of conformity of the transition, we divide an open and dense subset of ``$(2-1)$ with complex eigenvalues" into two open subsets. On one of those subsets, there is no topological invariants. On the other, we prove the existence of two moduli of stability. 
 
We prove here that a similar phenomenon occurs with only one complex eigenvalue.  

Before stating precisely our results, we should describe precisely the open and dense subset of $\cC$ we consider. Mainly, we assume that the saddles we consider are $C^1$-linearisable and we make some transversality assumptions.  

\subsection{Description of the set of connections}
We consider $p$ and $q$, two hyperbolic saddles in dimension three and we assume that for each saddle there exists a neighbourhood and $C^1$-coordinates (on this neighbourhood) such that the vector field is linear. 

\begin{remark}
In the set of hyperbolic saddles, the subset of $C^1$-linearisable saddles is open and dense (see for example \cite{be,be2,sa}). In dimension three, every saddle with complex eigenvalues is $C^1$-linearisable.
\end{remark}





Let us introduce some notations. We divide the set $(1-1)$ into $(1-1)_\CC$, $(1-1)_\RR$ according to the nature of the unstable eigenvalue of $q$ (complex meaning ``non-real"). Similarly, the set $(2-2)$ is the union of $(2-2)_\CC$ and $(2-2)_\RR$ according to the stable eigenvalue of $q$.

The set $(2-1)$ is divided into $(2-1)_{\CC,\CC}$, $(2-1)_{\RR,\CC}$, $(2-1)_{\CC,\RR}$ and $(2-1)_{\RR,\RR}$ with respect to the nature of the stable eigenvalues of $p$ and the unstable eigenvalues of $q$ respectively.

Besides the assumptions of hyperbolicity and linearization of the saddles, we will restrict our study of $\cC$ to the following generic subset: 

\begin{enumerate}
\item \label{a.first} For all the connections, we assume that a saddle with real eigenvalues in the stable direction at $p$ or in the unstable direction at $q$ has no pair of equal eigenvalues: there exist a strong and a weak stable (resp. unstable) direction in $W^s_{loc}(p)$ (resp. $W^u_{loc}(q)$) if the eigenvalues are real.

\item We define ${(1-2)}^g$ as the subset of $(1-2)$ such that the intersection of the unstable manifold of $p$ and the stable manifold of $q$ is transverse. 

\item For a connection in $(2-2)_\RR$, the stable eigenvalues of $p$ are real, and different by assumption~(\ref{a.first}). Let us denote by $E_p^{s,u}$  the 2-plane of $T_pM$ spanned by the weak stable direction and the unstable direction in $p$. We denote by $E_p^{ss,u}$ the 2-plane of $T_pM$ spanned by the strong stable direction and the unstable direction in $p$. The linearization assumption implies that the vector field is tangent to those planes and their intersection is included in the connection. 

We assume that the stable manifold of $q$ is transverse to $E_p^{s,u}$
and $E_p^{ss,u}$. We denote by $(2-2)^g_\RR$ the subset of $(2-2)_\RR$ with connections satisfying this assumption.  

\item Similarly, for a connection in $(1-1)_\RR$, we define $F_q^{u,s}$ and $F_q^{uu,s}$. 

We assume that the unstable manifold of $p$ is transverse to $F_q^{u,s}$ and $F_q^{uu,s}$. The subset obtained is denoted by $(1-1)^g_\RR$.

\item \label{a.cc} We assume that for connections in $(2-1)_{\CC,\CC}$ the ratio $\frac{-\Re(w^s(p))}{\Im(w^s(p))}$ is different from the ratio $\frac{\Re(w^u(q))}{\Im(w^u(q))}$, where $w^s(p)$ and $w^u(q)$ are the stable and unstable eigenvalues at $p$ and $q$ respectively.

\item \label{a.rr} \label{a.last} The last assumption about connections in $(2-1)$ concerns $(2-1)^g_{\RR,\RR}$, subset of $(2-1)_{\RR,\RR}$.

Consider $x$ be a point of $\Gamma\setminus\{p,q\}$. There exists a unique plane $E_x^{s,u}$ (resp. $F_x^{u,s}$) in $T_xM$ such that $E_x^{s,u}=T_x\{\phi^t(E_p^{s,u}),\,t\in \RR^+\}$ (resp. $F_x^{u,s}=T_x\{\phi^t(F_q^{u,s}),\,t\in \RR^-\}$).
Similarly, there exist a unique plane $E_x^{ss,u}$ (resp. $F_x^{uu,s}$) in $T_xM$ such that $E_x^{ss,u}=T_x\{\phi^t(E_p^{ss,u}),\,t\in \RR^+\}$ (resp. $F_x^{uu,s}=T_x\{\phi^t(F_q^{uu,s}),\,t\in \RR^-\}$).

We assume that for any $x$ in $\Gamma\setminus\{p,q\}$, $E_x^{s,u}$ and $F_x^{u,s}$ intersect transversally in $T_xM$. We assume also that $E_x^{ss,u}$ and $F_x^{u,s}$ (resp. $E_x^{s,u}$ and $F_x^{uu,s}$) intersect transversally. 
\end{enumerate}

\begin{remark}
The assumptions of transversality for any $x$ in $\Gamma\setminus\{p,q\}$ hold as soon as they hold for one $x_0$ in $\Gamma\setminus\{p,q\}$.
\end{remark}

We denote by $(2-1)^g$ the union of $(2-1)^g_{\CC,\CC}$, $(2-1)^g_{\RR,\CC}$, $(2-1)^g_{\CC,\RR}$ and ${(2-1)^g_{\RR,\RR}}$ and by $\cC^g$ the union of ${(1-2)}^g$, $(1-1)^g$, $(2-2)^g$ and $(2-1)^g$. The subset $\cC^g$ is an intersection of open and dense subsets of $\cC$.

\begin{remark}
With the assumptions made above, connections in ${(1-2)}^g$ are of  codimension $0$, connections in $(1-1)^g$ and $(2-2)^g$ are of codimension $1$ and connections in $(2-1)^g$ are of codimension $2$.
\end{remark}

\subsection{Description of connections in $(2-1)^g$} 
Let us denote by $U_p$ and $U_q$ two neighbourhoods of $p$ and $q$ respectively where the vector field is linear. We consider $\Sigma_p$ and $\Sigma_q$ two disks in $U_p$ and $U_q$, transverse to $\Gamma$ at $x_p$ and $x_q$. The Poincar\'e map $\cT\colon \Sigma_q \to\Sigma_p$ is a diffeomorphism at $x_p$ (in order to keep the same notations as in~\cite{bodu}, $\cT$ is in fact the inverse of the natural Poincar\'e map). 

\begin{lemma}\label{l.coord}
There exist $\Sigma_p$, $\Sigma_q$ and $C^1-$linearising coordinates in $U_p$ and $U_q$ such that $T$, the differential of $\cT$ at $x_p$, is given by $T=B_0\circ B_\lambda\circ B_1$ where $B_0$ and $B_1$ are rotations in $\Sigma_p$ and $\Sigma_q$ respectively and $B_\lambda$ is the matrix:

$$B_\lambda= \left(\begin{array}{cc}1&0\\0&\lambda\end{array}\right), \quad \lambda\geq1$$ (in the coordinates restricted to $\Sigma_p$ and $\Sigma_q$).
\end{lemma}

Related to $\lambda$, the quantities we consider are $t=\frac 1 2 (\lambda +\frac 1 \lambda) \mbox{ et } s=\frac 1 2 (\lambda -\frac 1 \lambda).$

In the next two paragraphs, we define a class, \itt$\subset (2-1)^g$, depending on the (normalised) eigenvalues and on $\lambda$ (or $t$ or $s$). We give a geometric interpretation of the distinction between connections in \itt~and others in section~\ref{s.2-1}.

\subsubsection{Connections in $(2-1)^g_{\CC,\CC}$}\label{ss.bodu}
Let $w^s(p)$ and $w^u(q)$ be the stable and unstable eigenvalues at $p$ and $q$ respectively. We denote by $\alpha=\frac{-\Re(w^s(p))}{\Im(w^s(p))}$ and by $\beta=\frac{\Re(w^u(q))}{\Im(w^u(q))}$; by assumption~(\ref{a.cc}), $\alpha\neq\beta$. 

The linearization coordinates are invariant by rotations centered on the connection. We can improve lemma~\ref{l.coord} to write $T=B_\lambda$. One can prove that up to topological equivalence, a connection in $(2-1)^g_{\CC,\CC}$ is defined by $(\alpha,\beta,\lambda)$.

We define the function $\psi\colon(\RR^*)^2\to \RR$ by 
$\psi(\alpha,\beta)= \frac{-\alpha\beta+|\alpha \beta| \sqrt{(\alpha^2+1)({\beta}^2+1)}}{\alpha^2 \beta^2}$. We note that $\psi(\alpha,\beta) > 1$ as soon as $\alpha\neq\beta$.

\begin{definition}
A connection $\Gamma$ in $(2-1)^g_{\CC,\CC}$ belongs to the class \itt~if $t\leq\psi(\alpha,\beta)$
\end{definition}

Moreover, let $U=\{{(\alpha,\beta,t)\in\RR^3}, {\alpha\neq 0},\ {\beta\neq 0,\ {\beta\neq\alpha}\mbox{ and }  t>\psi(\alpha,\beta)}\}$, we construct, in section~\ref{ss.loglog}, a continuous function $\Psi\colon U\to\RR$ ($\Psi$ is therefore defined on the complement of the class \itt).

\subsubsection{Connections in $(2-1)^g_{\RR,\CC}$}
Let $\Gamma$ be a connection of $(2-1)^g_{\RR,\CC}$. We denote by $w^u(q)$ the unstable eigenvalues at $q$ and by $\mu^s$ and $\mu^{ss}$ the stable eigenvalues at $p$. We denote by $\beta=\frac{\Re(w^u(q))}{\Im(w^u(q))}$ and by $\mu=\frac{\mu^s}{\mu^{ss}}$.

Only the coordinates on $U_q$ are invariant by rotations. So up to topological equivalence, a connection in $(2-1)^g_{\RR,\CC}$ is defined by $(\beta,\mu,\lambda,\theta_0)$, $\theta_0$ being the angle of the rotation $B_0$ of lemma~\ref{l.coord}.

Let us define $\mu_\pm=1+2\beta^2\pm 2|\beta|\sqrt{1+\beta^2}$ and $s_{\pm}=\frac{-\beta(\mu+1)\pm 2 |\beta|\sqrt{\mu(\beta^2+1)}}
{(\mu-1)\cos\theta_0\sin\theta_0}$. We denote $s_m$ the smaller and $s_M$ the larger of $s_+$ and $s_-$. We recall that $s=\frac 1 2 (\lambda -\frac 1 \lambda)$.

\begin{definition}\label{d.rctt}
A connection $\Gamma$ of $(2-1)^g_{\RR,\CC}$ belongs to the class \itt~if 
\begin{itemize}
\item $\theta_0=k \frac \pi 2$ ($k$ in $\ZZ$) and $\mu \in [\mu_-, \mu_+]$;
	
\item $\theta_0\neq k \frac \pi 2$ (for every $k$ in $\ZZ$) and 
\begin{itemize}
\item either $\mu \in ]\mu_-, \mu_+[$ and $s\leq s_M$,
\item or $\mu \notin [\mu_-, \mu_+]$ and $s_m$ is positive and $s \in [s_m,s_M]$,	
\item or $\mu=\mu_+$ or $\mu=\mu_-$ and $s\leq \frac{-2(\mu+1)}{\beta \cos \theta_0 \sin\theta_0 (\mu -1)}$. 
\end{itemize}
\end{itemize}
\end{definition}

On the complement of the class \itt~in $(2-1)^g_{\RR,\CC}$, we define a function $\Upsilon(\beta,\mu,s,\theta_0)$, the same way we define $\Psi$ for $(2-1)^g_{\CC,\CC}$ (see section~\ref{ss.upsilon}).

\subsubsection{Connections in $(2-1)^g_{\CC,\RR}$}
Let $\Gamma$ be a connection of $(2-1)^g_{\RR,\RR}$. We denote by $w^s(p)$ the stable eigenvalues at $p$ and by $\gamma^u$ and $\gamma^{uu}$ the unstable eigenvalues at $q$. We denote by $\alpha=\frac{\Re(w^s(p))}{\Im(w^s(p))}$ and by $\gamma=\frac{\gamma^s}{\gamma^{ss}}$.

If $X$ denotes the vector field in a neighbourhood of $\Gamma$, the vector field $-X$ presents a connection $-\Gamma$ in $(2-1)^g_{\RR,\CC}$. The parameters of $-\Gamma$ with the notation of the previous paragraph are $(\beta=\alpha,\mu=\gamma,\frac 1\lambda,\theta_0=-\theta_1)$.

\begin{definition}
A connection $\Gamma$ in $(2-1)^g_{\CC,\RR}$ belongs to the class $(TT)$ if the connection $-\Gamma$ in $(2-1)^g_{\RR,\CC}$ does.
\end{definition}

\subsubsection{Connections in $(2-1)^g_{\RR,\RR}$}
Let $x$ be a point of $\Gamma\setminus\{p,q\}$, for $\Gamma$ in $(2-1)^g_{\RR,\RR}$. And  let $\Sigma$ be a disk transverse to $\Gamma$ at $x$. The planes $E_x^{s,u}$, $F_x^{u,s}$, $E_x^{ss,u}$ and $F_x^{uu,s}$  give, by intersection with $T_x\Sigma$, four points in $\RR P(1)\approx S^1$, denoted by $\omega_p^s$, $\omega_q^u$, $\omega_p^{ss}$ and $\omega_q^{uu}$ respectively.

By assumption~(\ref{a.rr}), we have $\omega_p^s\neq\omega_q^u$, and the points $\omega_p^{ss}$, $\omega_q^{uu}$ belong to $\RR P(1)\setminus\{\omega_p^s,\omega_q^u \}$.

\begin{definition}
A connection $\Gamma$ in $(2-1)^g_{\RR,\RR}$ is of type $(I)$ if $\omega_p^{ss}$ and $\omega_q^{uu}$ are in the same connected component of $\RR P(1)\setminus\{\omega_p^s,\omega_q^u \}$. Otherwise, it is of type $(II)$.
\end{definition}


\subsection{Statement of the results}

Up to topological equivalence, the Morse index of a saddle is a complete invariant in a neighbourhood of the saddle (Theorem of Hartman-Grobman). Therefore, if two connections belong to two different subsets $(1-2)^g$, $(1-1)^g$, $(2-2)^g$ and $(2-1)^g$, they are not topologically equivalent.
If they are in the same subset, we have the following three results:

\begin{theorem}\label{t.1-2}
Two connections of $(1-2)^g$ are always topologically equivalent.
\end{theorem}

\begin{theorem}\label{t.1-1}
Two connections of $(1-1)^g$ are topologically equivalent if and only if they are in the same subset $(1-1)^g_\RR$ or $(1-1)^g_\CC$.

Two connections of $(2-2)^g$ are topologically equivalent if and only if they are in the same subset $(2-2)^g_\RR$ or $(2-2)^g_\CC$.
\end{theorem}

\newpage

\begin{theorem}\label{t.2-1}

Two connections $\Gamma$ and $\tilde \Gamma$ of $(2-1)^g$ are topologically equivalent if and only if

$$\hspace{-.2cm}
\begin{array}{|c||c|c|c|c|}\hline
\qquad \Gamma \in &&&&\\ 
      &  (2-1)^g_{\CC,\CC}& (2-1)^g_{\RR,\CC}&(2-1)^g_{\CC,\RR}&(2-1)^g_{\RR,\RR} \\
\tilde \Gamma \in \qquad&&&&\\     
\hline
\hline
&\mbox{\scriptsize $\Gamma$ \& $\tilde \Gamma \in$ \itt} \mbox{ or} &\mbox{\footnotesize $\Gamma \in$ \itt}&\mbox{\footnotesize $\Gamma \in$ \itt}&\\ 

(2-1)^g_{\CC,\CC}&\mbox{\scriptsize \big($\Gamma${\footnotesize \& }$\tilde \Gamma \notin$ \itt~ {\footnotesize \& } $\frac\alpha\beta=\frac{\tilde\alpha}{\tilde\beta}$} 
&\mbox{\&} &\mbox{\&}& \mbox{\footnotesize impossible} \\

 &\mbox{\footnotesize{\& $\Psi(\alpha,\beta,t)=\Psi(\tilde \alpha,\tilde \beta,\tilde t)$\big)}} &\mbox{\footnotesize $\tilde \Gamma \in$ \itt}  &\mbox{\footnotesize $\tilde \Gamma \in$ \itt} & \\

\hline 
&&&&\\
&\mbox{\footnotesize $\Gamma \in$ \itt}&\mbox{\scriptsize $\Gamma$ \& $\tilde \Gamma \in$ \itt} \mbox{  or} &\mbox{\footnotesize $\Gamma \in$ \itt}&\\

(2-1)^g_{\RR,\CC}&\mbox{\&}&\mbox{\scriptsize \big($\Gamma${\footnotesize \& }$\tilde \Gamma \notin$ \itt~ {\footnotesize \& }} &\mbox{\&} & \mbox{\footnotesize impossible} \\ 

& \mbox{\footnotesize $\tilde \Gamma \in$ \itt}&\mbox{\scriptsize{$\Upsilon(\beta,\mu,s,\theta_0)=$}}&\mbox{\footnotesize $\tilde \Gamma \in$ \itt} & \\
&&\mbox{\scriptsize{$\Upsilon(\tilde \beta,\tilde \mu,\tilde s,\tilde \theta_0)$}\big)} &&\\
\hline
&&&&\\
&\mbox{\footnotesize $\Gamma \in$ \itt }&\mbox{\footnotesize $\Gamma \in$ \itt}&\mbox{\scriptsize $\Gamma$ \& $\tilde \Gamma \in$ \itt} \mbox{  or}&\\

 (2-1)^g_{\CC,\RR}&\mbox{\&} &\mbox{\&}&\mbox{\scriptsize \big($\Gamma${\footnotesize \& }$\tilde \Gamma \notin$ \itt~ {\footnotesize \& }}& \mbox{\footnotesize impossible} \\ 

& \mbox{\footnotesize $\tilde \Gamma \in$ \itt}&\mbox{\footnotesize $\tilde \Gamma \in$ \itt}  &\mbox{\scriptsize{$\Upsilon(\alpha,\gamma,-s,-\theta_1)=$}}& \\  

& & & \mbox{\scriptsize{$\Upsilon(\tilde \alpha,\tilde \gamma,-\tilde s,-\tilde \theta_1)$}\big)}& \\
\hline 
&&&& \mbox{\footnotesize $\Gamma$ \& $\tilde \Gamma$ have}\\
(2-1)^g_{\RR,\RR}&\mbox{\footnotesize impossible} &\mbox{\footnotesize impossible}&\mbox{\footnotesize impossible}& \mbox{\footnotesize same type}\\
& & & & \mbox{\footnotesize ((I)\mbox{ or } (II))}\\
\hline	
\end{array}
$$
\label{t.principal}
\end{theorem}

We proved in~\cite{bodu} that for a connection in $(2-1)^g_{\CC,\CC}$, the case $\alpha=\beta$ leads to another class of topological equivalence.

Let us mention the work of Clementa Alonso~\cite{al} where the case $(2-1)^g_{\RR,\RR}$ with $E_x^{s,u}=F_x^{u,s}$ and $E_x^{ss,u}=F_x^{uu,s}$ is completely classified (for analytic vector fields), using different methods.

Finally, we remark that \cite[Theorem 2.1]{ducras1} implies that the absolute value of the imaginary part of the complex eigenvalue of a saddle of a connection in $(1-1)^g_\CC$ or $(2-2)^g_\CC$ is a conjugacy invariant.

The organisation of the paper is as follow. First we investigate the behaviour around one saddle (section~\ref{s.uneselle}), in particular we give necessary and sufficient conditions for a homeomorphism between two sections to be extended into a topological equivalence (sections~\ref{ss.ind2} and \ref{ss.ind1}, compare~\cite{al,bodu}). We defined in~\cite{bodu} an equivalence relation (to be {\em semblable}) on the set of radial-type foliations and of vertical foliations of the disk. And we associate to a saddle an equivalence class of foliations corresponding to the equivalence class up to topological equivalence of the saddle (section~\ref{s.connect}). Then we prove the existence of {\em characteristic foliations} (section~\ref{ss.charac}) for every connections in $\cC^g$ such that topological equivalence in a neighbourhood of the connection is equivalent to the equivalence (up to the semblable relation) of a pair of characteristic foliations. Then our problem is reduced to the classification of pairs of characteristic foliations, this is achieved in section~\ref{s.class}.

\quad

{\bf Acknowledgements:} This work is part of my thesis dissertation, completed under the guidance of Robert Roussarie. I would like to thank Patrice Le Calvez for many improvements he suggested and Floris Takens who encouraged me to write those results in this form. I would also like to thank C. Alonso for many discussions on this subject. Finally, I would like to thank Sebastian van Strien who gave me very useful comments and made possible my stay at Warwick University.

\section{One saddle, blow-up}\label{s.uneselle}
In this section, we investigate the topological behaviour of a vector field in a neighbourhood of a hyperbolic saddle. In particular, we give sufficient and necessary conditions to extend a homeomorphism into a topological equivalence.

\subsection{Blow-up of a radial foliation on a disk}

A foliation $\cF$ of a disk $D$, singular at a point $A$, is a {\em foliation of radial-type} of the disk $(D,A)$ if $\cF$ is homeomorphic to the linear radial foliation of $D$ in A. We will assume that $\cF$ is transverse to $\partial D$.

\begin{definition}
A {\em blow-up} of $\cF$ is a homeomorphism $f \colon D\setminus\{A\} \to S^1\times ]0,1]$, such that the foliation $f(\cF)$ extends to a foliation of $S^1\times [0,1]$, homeomorphic to the vertical foliation of the annulus.
\end{definition}



Consider $\cF$ and $\tilde \cF$ two foliations of radial-type of $(D,A)$ and, $f$ and $\tilde f$, blow-up of $\cF$ and $\tilde \cF$ respectively.

\begin{definition}\label{d.guide1}
The foliation $\cF$ {\em guide} $\tilde \cF$ if $\tilde f\circ f^{-1}$ extends to a continuous map of $S^1\times [0,1]$.

The foliations $\cF$ and $\tilde \cF$ are {\em semblable} if $\cF$  guide $\tilde \cF$ and $\tilde \cF$ guide $\cF$. This will be denoted by $\cF \simeq \tilde\cF$.
\end{definition}

The relation ``to be semblable" is an equivalence relation. In order to read the notion of guidance directly on the disk, we define a {\em sector} of a radial foliation $\cF$ as the union of leaves cutting an interval $I$ in $\partial D$. A sector is {\em open} (resp. {\em closed}) if the interval $I$ is open (resp. closed) in $\partial D$. Definition~\ref{d.guide1} is equivalent to the following (compare~\cite{bodu}).

\begin{definition}
The foliation $\cF$ {\em guide} $\tilde \cF$ if for every covering $\cU$ of $D\setminus{A}$ by open sectors of $\cF$, there exists a neighbourhood $V$ of $A$ such that for every leaf $\tilde L$ of $\tilde F$, there exists a sector $S$ in $\cU$ with $\tilde L \cap V \subset S$.

The foliations $\cF$ and $\tilde \cF$ are {\em semblable} if $\cF$  guide $\tilde \cF$ and $\tilde \cF$ guide $\cF$.
\end{definition}

It is easy to see that the definition of {\em guide} does not depend on the blow-up and that we have the following.
 
\begin{lemma}
The foliations $\cF$ and $\tilde \cF$ are semblable if and only if $\tilde f\circ f^{-1}$ extends in a homeomorphism of $S^1\times [0,1]$ (i.e. $\cF$ and $\tilde \cF$ have the same blow-up).
\end{lemma}

Let $\cF$ (resp. $\tilde \cF$) be a radial-type foliation of the disk $(D,A)$ (resp. $(\tilde D,\tilde A)$), let $f$ (resp. $\tilde f$) be a blow-up of $\cF$ (resp. $\tilde \cF$) and consider a homeomorphism $h\colon D \to \tilde D$ such that $h(A)=\tilde A$. 

\begin{lemma}\label{l.blow}
The map $\tilde f\circ h$ is a blow-up of $\cF$ if and only if $h(\cF)$ and $\tilde \cF$ are semblable.
\end{lemma}

In the next two sections, we give necessary and sufficient conditions for a homeomorphism between sections of vector fields in a neighbourhood of a saddle to be extended into a topological equivalence.

\subsection{Adapted foliation of a saddle of Morse index 2}\label{ss.ind2}
Let $p$ be a saddle point of Morse index $2$ of a vector field in dimension three. Let $D$ be a disc transverse to the flow at $A$, a point of the local unstable manifold of $p$. We consider $C$ an annulus (diffeomorphic to $S^1\times[0,1]$), transverse to the flow, such that $S^1\times\{0\}$ is included in the local stable manifold of $p$.

There exists a neighbourhood $U$ of $S^1\times\{0\}$ in $C$ and a neighbourhood $V$ of $A$ in $D$ such that the Poincar\'e map $P$ between $U\setminus S^1\times\{0\}$ and $V\setminus \{A\}$ is well-defined. 

\begin{definition}
The image under $P$ of the vertical foliation of $C$ gives a radial-type foliation of $(D,A)$ called an {\em adapted foliation} for $p$.
\end{definition}

From last section, it is clear that different choices of disks $D$, or annuli $C$ would lead to semblable adapted foliations. From lemma~\ref{l.blow} we deduce the following.

\begin{lemma}\label{l.equivalence}
Let $\cF$ and $\tilde \cF$ be adapted foliations of $p$ and $\tilde p$ on disks $D$ and $\tilde D$ respectively.

A homeomorphism $h\colon D\to \tilde D$ can be extended into a topological equivalence between neighbourhoods of the saddles if and only if $h(\cF)$ is semblable to $\tilde \cF$.
\end{lemma}

\subsection{Adapted foliation of a saddle of Morse index 1}\label{ss.ind1}
Let $p$ be a saddle point of a vector field in dimension three of Morse index $1$. We consider $C$ an annulus, transverse to the flow, such that $S^1\times\{0\}$ is included in the local unstable manifold of $p$.



\begin{lemma}
Consider $p$, $\tilde p$ two saddles of Morse index one and $p_1$, $\tilde p_1$ two points of $C\cap W^u_{loc}(p)$ and $\tilde C \cap W^u_{loc}(\tilde p)$ respectively. Let $N$ and $\tilde N$ be two disks, neighbourhood in $C$ and $\tilde C$ of $p_1$ and $\tilde p_1$ respectively.
Every homeomorphism $h\colon N \to\tilde N$ mapping $N\cap W^u_{loc}(p)$ to $\tilde N \cap W^u_{loc}(\tilde p)$ can be extended into a topological equivalence.
\end{lemma}

\begin{proof}
One can extend $h$ into a homeomorphism from a neighbourhood of $C\cap W^u_{loc}(p)$ in $C$ to a neighbourhood of $\tilde C \cap W^u_{loc}(\tilde p)$ in $\tilde C$. We can extend $h$ such that it maps $C\cap W^u_{loc}(p)$ to $\tilde C \cap W^u_{loc}(\tilde p)$ and the vertical foliation of a a neighbourhood of $\tilde C \cap W^u_{loc}(\tilde p)$ to a vertical foliation in $\tilde C$.

Let us consider $\Sigma$ (resp. $\tilde \Sigma$) a disk transverse to the local stable manifold of $p$ (resp. local stable manifold of $\tilde p$) at a point $x\neq p$ (resp. $\tilde x\neq \tilde p$).
Using the Poincar\'e maps $P^{-1}\colon C\to \Sigma$ and $\tilde P^{-1}\colon \tilde C\to \tilde\Sigma$, we can define $h_0$ from $\Sigma$ to $\tilde \Sigma$ from $h$. The image by $P^{-1}$ (resp. $\tilde P^{-1}$) of the vertical foliation of a neighbourhood of $C \cap W^u_{loc}(p)$ (resp. $\tilde C \cap W^u_{loc}(\tilde p)$) is an adapted foliation for the saddle $p$ (resp. $\tilde p$) with reversed time. And by construction, $h_0\colon\Sigma\to\tilde \Sigma$ maps one foliation to the other. By lemma~\ref{l.equivalence}, one can extend $h_0$ in a topological equivalence between $p$ and $\tilde p$ with reversed time so one can extend $h\colon N \to\tilde N$ in a topological equivalence.
\end{proof}

Let $p_1$ be a point of $C\cap W^u_{loc}(p)$ and $N$ be disk, neighbourhood of $p_1$ in $C$. To be coherent with the case of index two, we define:

\begin{definition}
An {\em adapted} foliation on $N$ is a vertical foliation based on $S=N\cap W^u_{loc}(p)$.
\end{definition}
\begin{definition}
Two adapted foliations (on $N$ and $\tilde N$) are {\em semblable} if there exists a homeomorphism from $N$ to $\tilde N$ mapping $S$ onto $\tilde S$.
\end{definition}

\section{Topological equivalence of saddle connections, characteristic foliations}\label{s.connect}
Let $\Gamma$ be a saddle connection in $\cC^g$. Let $A$ be a point of $\Gamma \setminus\{p,q\}$ and $D$ be a disk transverse to the flow in $A$. For a saddle $p$ of a vector field $X$, we denote by $-p$ the saddle of the vector field $-X$. On $D$, we consider $\cF_p$ an adapted foliation for $-p$ and $\cF_q$ an adapted foliation for $q$.

\begin{lemma}
Two connections $\Gamma$ and $\tilde \Gamma$ are topologically equivalent if and only if there exists $D$ and $\tilde D$ as above and $h\colon D \to \tilde D$ such that $h(\cF_p)$ is semblable to $\tilde \cF_p$ and $h(\cF_q)$ is semblable to $\tilde \cF_q$. 
\label{l.disque}
\end{lemma}

The classification up to topological equivalence of the connections in $\cC^g$ is equivalent to the classification of pairs of adapted foliations on a disk with respect to the ``semblable" relation. In order to achieve the later classification, we choose ``good" adapted foliations: the characteristics foliations.

\subsection{Characteristics foliations}\label{ss.charac}

\begin{definition}
A {\em vertical} foliation of the disk with the segment $S$ as a {\em basis} is a $C^1$-foliation given by a tubular neighbourhood's vertical foliation of $S$ in the disk. The image by a diffeomorphism of a vertical foliation is still named vertical foliation. 
\end{definition}

On the plane $(x,y)$, let $R(x,y)=x\frac\partial{\partial x}+y\frac\partial{\partial y}$ denotes the radial vector field and $\frac\partial{\partial \theta}=-y\frac\partial{\partial x}+x\frac\partial{\partial y}$.
\begin{definition}
A {\em logarithmic} foliation $\cL_{\alpha}$ of the disk is a foliation obtained by integrating the vector field: $X_\alpha=R+\alpha\frac\partial{\partial \theta}$.
\end{definition}

\begin{definition}
A {\em real} foliation $\cR_{\gamma}$ is given by the vector field: $Y_\gamma=x\frac\partial{\partial x}+\gamma y\frac\partial{\partial y}$.
The image by a diffeomorphism of a real foliation is still named real foliation.
\end{definition}

For logarithmic foliations, we have the following proposition~:

\begin{proposition}[\cite{bodu}, Proposition 2.1]
Let $\cL_\alpha$ be a logarithmic foliation of $\RR^2$ and let $\varphi\colon\RR^2\to\RR^2$ be a diffeomorphism, fixing the origin $(0,0)$ and tangent to the identity at $(0,0)$. Then we have $\varphi(\cL_\alpha)\simeq\cL_\alpha$.
\end{proposition}

We make use of the assumptions made in the introduction (existence of $C^1-$linearising coordinates for saddles in $\cC^g$ and assumptions (\ref{a.first})-(\ref{a.last})) to obtain the following lemma. We recall that $B_\lambda$ denotes the matrix $B_\lambda=\left(\begin{array}{cc} 1&0\\ 0&\lambda\end{array}\right), \mbox{ with } \lambda\geq 1$.

\begin{lemma}
Let $\Gamma$ be a connection in $\cC^g$. There exist a disk $D$, transverse to $\Gamma$ in $A$, and $C^1$ coordinates on $D$ such that the following foliations are adapted to $-p$ and $q$ respectively:
\begin{itemize}
\item {$\Gamma \in (1-2)^g$.} Two vertical foliations with basis $S_p$ and $S_q$ intersecting transversally in $A$.
\item {$\Gamma \in (1-1)^g_\CC$:} one vertical foliation and one logarithmic foliation $\cL_\beta$. 
\item {$\Gamma \in (1-1)^g_\RR$:} one vertical foliation and one real foliation $\cR_\mu$ such that the basis of the vertical foliation is transverse to the weak and strong manifolds of the real foliation.
\item {$\Gamma \in (2-2)^g_\CC$:} one logarithmic foliation $\cL_\alpha$ and one vertical foliation.
\item {$\Gamma \in (2-2)^g_\RR$:} one real foliation $\cR_\gamma$ and one vertical foliation such that the basis of the vertical foliation is transverse to the weak and strong manifolds of the real foliation.
\item {$\Gamma \in (2-1)^g_{\CC,\CC}$:} one logarithmic foliation $\cL_\alpha$ and the image by $B_\lambda$ of a logarithmic foliation $\cL_\beta$. We denote by $(\cL_\alpha,\cL_\beta^\lambda)$ this pair of foliations. 
\item {$\Gamma \in (2-1)^g_{\RR,\CC}$:} one real foliation $\cR_\mu$ and the image by a rotation of angle $\theta_0$ of the image by $B_\lambda$ of a logarithmic foliation $\cL_\beta$. This is denoted by $(\cR_\mu,\cL_\beta^{\lambda,\theta_0}={B_0}^\star({B_\lambda}^\star(\cL_\beta)))$.
\item {$\Gamma \in (2-1)^g_{\CC,\RR}$:} with the notations above, $(\cL_\alpha^{\frac 1\lambda,-\theta_1},\cR_\gamma)$.

\item {$\Gamma \in (2-1)^g_{\RR,\RR}$:} one real foliation $\cR_\mu$ and the image by a diffeomorphism of $\cR_\gamma$, such that the weak invariant manifolds are transverse, and the strong invariant manifolds are transverse to the weak ones.
\end{itemize}
\label{l.liste}
\end{lemma}

We call {\em complex} foliation the image by $B_\lambda$ and a rotation of a logarithmic foliation.

\begin{definition}
We call {\em characteristic} foliations of $\Gamma$ the pair of foliations assigned by lemma~\ref{l.liste} to the saddles of the connection.
\end{definition}

\begin{definition}
Two pairs of foliations $(\cF,\cG)$ and $(\tilde\cF,\tilde\cG)$ are {\em semblable} if there exists a homeomorphism $h$ such that $h(\cF)$ is semblable to $\tilde\cF$ and $h(\cG)$ is semblable to $\tilde\cG$.
\end{definition}

Combining lemma~\ref{l.disque} and lemma~\ref{l.liste}, we obtain:
\begin{lemma}\label{l.charac}
Two connections $\Gamma$ and $\tilde \Gamma$ are topologically equivalent if and only if the characteristic foliations of $\Gamma$ are semblable to the characteristic foliations of $\tilde \Gamma$.
\end{lemma}

\section{Classification of semblable pairs of characteristic foliations}\label{s.class}
We consider two pairs of characteristic foliations and we give necessary and sufficient conditions for the existence of a homeomorphism mapping the first pair to a pair that is semblable to the second pair of foliations.

\subsection{Two vertical foliations}

Let us consider $\cV_p$ and $\cV_q$ two vertical foliations with basis $S_p$ and $S_q$ respectively. We assume that $S_p$ and $S_q$ intersect transversally in $A$. Next lemma, together with lemma~\ref{l.charac} and lemma~\ref{l.liste}, gives the proof of the Theorem~\ref{t.1-2}.

\begin{lemma}
Given two pairs of vertical foliations $(\cV_p,\cV_q)$ and $(\tilde \cV_p,\tilde \cV_q)$ as above there exist $U$ and $\tilde U$, neighbourhoods of $A$ and $\tilde A$ respectively and $h\colon U\to \tilde U$ such that $h(\cV_p)=\tilde \cV_p$ and $h(\cV_q)=\tilde \cV_q$.
\end{lemma}
\begin{proof}
There exists $U_0$ a neighbourhood of $A$ such that $\cV_p$ and $\cV_q$ are transverse in $U_0$. The same holds for $\tilde U_0$, neighbourhood of $\tilde A$. Define $h\colon (S_p\cup S_q)\cap U_0\to (\tilde S_p\cup \tilde S_q)\cap \tilde U_0$ a arbitrary homeomorphism $(h(A)=\tilde A)$.
The foliations $\cV_p$ and $\cV_q$ give continuous projections $\pi_p$ and $\pi_q$ from $U_0$ to $S_p$ and $S_q$. There exists $U\subset U_0$ such that we extend $h$ to  by $h(x)=\pi_{\tilde p}^{-1}(h(\pi_p(x)))\cap\pi_{\tilde q}^{-1}(h(\pi_q(x)))$ into a homeomorphism onto its image $\tilde U$.
\end{proof}

\subsection{One vertical foliation}
Let us prove Theorem~\ref{t.1-1} in two steps. Let $(\cV,\cF)$ be a pair of foliations such that $\cV$ is a vertical foliation based on $S$. The {\em nature} of the pair $(\cV,\cF)$ is real if $\cF$ is a real foliation and complex if $\cF$ is a complex foliation.

\begin{lemma}
For two pairs of foliations $(\cV,\cF)$ and $(\tilde \cV,\tilde \cF)$ with the same nature, there exists $U$ a neighbourhood of $A$, $\tilde U$ a neighbourhood of $\tilde A$ and a homeomorphism $h\colon U\to \tilde U$ mapping $\cV$ on a foliation semblable to $\tilde V$ and $\cF$ on $\tilde \cF$.
\end{lemma} 
\begin{proof}
First, let us prove the lemma with complex foliations. Let $U$ be a neighbourhood of $A$ such that every leaf of $\cF$ is transverse to $S$. Consider $F$ one leaf of $\cF$, we orient $F$ ``to the origin". Let $F_0$ be a segment of $F$ starting on $x_0\in S\cap U$ and finishing on $x_2\in S$ the second intersection of $F$ with $S$ after $x_0$. Similarly, we consider an interval $[\tilde x_0, \tilde x_2]$ on $\tilde S$. Let $h\colon[x_0,x_2]\to[\tilde x_0, \tilde x_2]$ a homeomorphism such that $h(x_0)=\tilde x_0$ and $h(x_1)=\tilde x_1$. We set $h(A)=\tilde A$. Using a parametrisation of the leaves of $\cF$ and $\tilde \cF$ we extend $h$ to $U$ such that $h(\cF)=\tilde \cF$ and $h(S)= \tilde S$ so $h(\cV)\simeq\tilde\cV$.

Second, assume that $\cF$ is a real foliation. Let $U$ (resp. $\tilde U$) be a neighbourhood of $A$ (resp. $\tilde A$), and $h$ be a homeomorphism from $S\cap U$ to $\tilde S\cap\tilde U$, such that $h(A)=\tilde A$. The invariant manifolds of $\cF$ divide the plane into two cones, we assume that in $U$, $S$ is included in one cone. We extend $h$ in this cone using a parametrisation of the leaves of $\cF$. On the other cone, we set $h$ to map $\cF$ on $\tilde \cF$ and to extend continuously to the invariant manifolds.
\end{proof}

\begin{lemma}\label{l.infini}
Two pairs of foliations $(\cV,\cF)$ and $(\tilde \cV,\tilde \cF)$ of different nature are not semblable.
\end{lemma}
\begin{proof}
Suppose there exist $(\cV,\cR)$ and $(\tilde \cV,\tilde \cL)$ such that $\cR$ is a real foliation and $\tilde \cL$ a complex one. And suppose there exists $h$ such that $h(\cV)\simeq\tilde \cV$ and $h(\cR)\simeq\tilde \cL$. In a neighbourhood of $A$, every leaf of $\cR$ cuts at most once $S$ whereas in a neighbourhood of $\tilde A$, every leaf of $\tilde \cL$ cuts $\tilde S$ infinitely many times. Let $\tilde L$ be a leaf of $\tilde \cL$ and $\tilde V$ a sector (for the foliation $\tilde \cL$) around $\tilde L$. There exists a leaf $R$ of $\cR$ such that $h(R)\cap \tilde U$ is included in $V$ for some $\tilde U$ neighbourhood of $\tilde A$. So $h(R)$ cuts infinitely many times $\tilde S$. As $h(S)=\tilde S$, this is a contradiction.
\end{proof}

The two previous lemmas proved the part of Theorem~\ref{t.1-1} concerning connections in $(1-1)^g$. The proof for connections in $(2-2)^g$ is similar.

\subsection{At least one complex foliation}\label{s.2-1}
Next definition explains the definition of the class \itt~for connections in $(2-1)^g$ (we remark that a pair of characteristic foliations with one complex foliation is invariant by homothety).

\begin{definition}
A pair of characteristic foliations with at least one complex foliation is of type \itt~if the foliations are everywhere (except the origin) topologically transverse. 
\end{definition}

For a pair of characteristic foliation, we can compute if its in \itt~or not.

\begin{lemma}
A pair of foliations $(\cL_\alpha,\cL_\beta^\lambda)$ ($\alpha\neq\beta\neq 0$) belongs to the class \itt~if and only if $t\leq \psi(\alpha,\beta)$.
\end{lemma}

\begin{proof}
The points of smooth tangencies are given by the equation: 
\begin{equation*} 
\label{e.toile} 
(x-\frac\beta\lambda y)(\alpha x +y)-(x-\alpha y)(\lambda \beta x + y)=0
\end{equation*} 

Foliations are invariant by homothety, so we can restrict to the line $\{ y= 1 \}$ to obtain: 
\begin{equation*}
(\alpha-\lambda \beta)x^2+(\alpha \beta \lambda - \frac{\alpha \beta}{\lambda})x+(\alpha-\frac{\beta}{\lambda})=0
\label{e.mmental}
\end{equation*}
Discriminant of the last equation with respect to $x$ is a degree two polynomial in $t$ (for $t \geq 1$): 
$$
\Delta(t)=4\big(\alpha^2 \beta^2 t^2 + 2 \alpha \beta t -(\alpha^2+\beta^2+\alpha^2 \beta^2)\big)
$$

To obtain the sign of $\Delta$, we have to solve $\Delta=0$ with respect $t$; we have a new discriminant: 
$$\Delta'=16\alpha^2 \beta^2(\alpha^2+1)(\beta^2+1)$$ 
 
We deduce the result: $\alpha\beta\neq 0$ then $\Delta' > 0$, so the equation $\Delta(t)=0$ has two roots in $\RR$. Only one root is greater than 1: $t_+ =    \frac{-\alpha\beta+|\alpha \beta| \sqrt{(\alpha^2+1)({\beta}^2+1)}}{\alpha^2 \beta^2}= \psi(\alpha,\beta).$ This allows to conclude.
\end{proof}

Let recall from the introduction that $\mu_\pm=1+2\beta^2\pm 2|\beta|\sqrt{1+\beta^2}$ and $s_{\pm}=\frac{-\beta(\mu+1)\pm 2 |\beta|\sqrt{\mu(\beta^2+1)}}
{(\mu-1)\cos\theta_0\sin\theta_0}$. We denote $s_m$ the smaller and $s_M$ the larger of $s_+$ and $s_-$, finally $s=\frac 1 2 (\lambda -\frac 1 \lambda)$.
\begin{lemma}
The pair of foliations $(\cR_\gamma,\cL_\alpha^{\lambda,\theta})$ belongs to\itt~if and only if

\begin{itemize}
\item $\theta_0=k \frac \pi 2$ ($k$ in $\ZZ$) and $\mu \in [\mu_-, \mu_+]$;
	
\item $\theta_0\neq k \frac \pi 2$ (for every $k$ in $\ZZ$) and 
\begin{itemize}
\item either $\mu \in ]\mu_-, \mu_+[$ and $s\leq s_M$,
\item or $\mu \notin [\mu_-, \mu_+]$ and $s_m$ is positive and $s \in [s_m,s_M]$,	
\item or $\mu=\mu_+$ or $\mu=\mu_-$ and $s\leq \frac{-2(\mu+1)}{\beta \cos \theta_0 \sin\theta_0 (\mu -1)}$. 
\end{itemize}
\end{itemize}
\end{lemma}

\begin{proof}
The proof is the same as the proof of previous lemma.
The points of tangencies are given by the following equation (with respect to $(x,y)$):
\begin{eqnarray*}
x\times ( \alpha(\frac {\sin^2\theta} \lambda+\lambda \cos^2\theta)x+(1+\alpha\sin\theta\cos\theta (\frac 1 \lambda-\lambda))y)-\qquad \qquad\qquad\nonumber\\
\qquad \qquad \qquad \gamma y \times ((1+\alpha\sin\theta\cos\theta(\lambda-\frac 1 \lambda))x-\alpha(\frac {\cos^2\theta} \lambda+\lambda \sin^2\theta)y)= 0. 
\label{e.reelcplx} 
\end{eqnarray*}
Again, we restrict to the horizontal $\{y=1\}$ to obtain the equation ix $x$:

\begin{equation*}
ax^2+bx+c=0,\label{e.deux}
\end{equation*}
with 
\begin{eqnarray*}
 a&=&\alpha (\frac{\sin^2\theta}\lambda+\lambda \cos^2\theta)\label{e.a},\\
 b&=&1+\alpha \cos\theta\sin\theta(\frac 1 \lambda - \lambda)-\gamma(1+\alpha \cos\theta\sin\theta(\lambda-\frac 1 \lambda))\ \label{e.b},\\
 c&=& \alpha \gamma (\frac{\cos^2\theta}\lambda+\lambda \sin^2\theta)\label{e.c}.
\end{eqnarray*}
With $s=\frac 1 2 (\lambda -\frac 1 \lambda)$, we can write down the discriminant as:
\begin{equation*}
\Delta=b^2-4ac=A {s}^2+Bs+C, 
\end{equation*}
with 
\begin{eqnarray*}
 A&=& \alpha^2\cos^2\theta \sin^2\theta(\gamma-1)^2\label{e.A},\\
 B&=& 2\alpha\sin\theta\cos\theta(\gamma^2-1)\label{e.B},\\
 C&=& \gamma^2-2\gamma(1+2\alpha^2)+1\label{e.C}.
\end{eqnarray*}
If $\Delta \leq0$ the foliations are topologically transverse on $\RR^2\setminus\{(0,0)\}$.
The discriminant of the last equation is given by:
$$\Delta'=B^2-4AC=16 \alpha^2 \cos^2\theta \sin^2\theta (\gamma-1)^2 \gamma(\alpha^2+1)$$

The real $\gamma$ is strictly positive, so $\Delta'\geq 0$. And $\Delta'=0$ if and only if $\theta=k \frac\pi 2$ with $k$ in $\ZZ$.

\begin{enumerate}
\item If $\theta=k \frac\pi 2$, $\Delta=C=\gamma^2-2\gamma(1+2\alpha^2)+1$, and $C(\gamma_\pm)=0$. Therefore, for every $\lambda$,  
\begin{enumerate}
\item if $\gamma \notin [\gamma_-,\gamma_+]$, $\Delta>0$;

\item if $\gamma \in [\gamma_-,\gamma_+]$, $\Delta\leq 0$. 
\end{enumerate}

\item If $\theta\neq k \frac\pi 2$, $\Delta'$ is strictly positive, the equation $\Delta=0$ has two roots, $$s_{\pm}=\frac{-\alpha(\gamma+1)\pm 2 \times 
|\alpha| \sqrt{\gamma(\alpha^2+1)}}{(\gamma-1)\cos\theta\sin\theta}.$$

The ratio $\frac C A$ has same sign as $C$. We have: 
\begin{enumerate}
\item if $\gamma \notin [\gamma_-,\gamma_+]$ then $C>0$. Then, either $s_\pm$ are negative roots of $\Delta=0$ and we have $s\geq 0$, $\Delta>0$. Or $s_\pm$ are positive roots of this equation and we have: if $s \notin [s_m,s_M]$, $\Delta>0$, if $s \in [s_m,s_M]$, $\Delta\leq 0$. 

\item If $\gamma \in ]\gamma_-,\gamma_+[$ then $C<0$, $s_M$ is positive and, if $s> s_M$, $\Delta >0$. If $s\leq s_M$ we have $\Delta \leq 0$.

\item Finally, if $\gamma=\gamma_\pm$ then $C=0$, and $s=0$ is a root of $\Delta=0$. In this case, if $\frac{-B} {A}=\frac{-2(\gamma+1)}{\alpha \cos \theta \sin\theta (\gamma -1)}$ is negative, $\Delta$ is positive for $s>0$ and $\Delta=0$ for $s=0$. Otherwise, $\Delta$ is negative if $s \leq \frac{-B} {A}$.
\end{enumerate}
\end{enumerate}
This ends the proof.
\end{proof}

We can summarise:
\begin{lemma}\label{l.calcul}
A connection in $(2-1)_{\CC,\CC}^g$, $(2-1)_{\RR,\CC}^g$ or $(2-1)_{\CC,\RR}^g$ belongs to the class \itt~if and only if its associated pair of characteristic foliations is of type \itt.
\end{lemma}

\begin{figure}[htb]
\centerline{\includegraphics{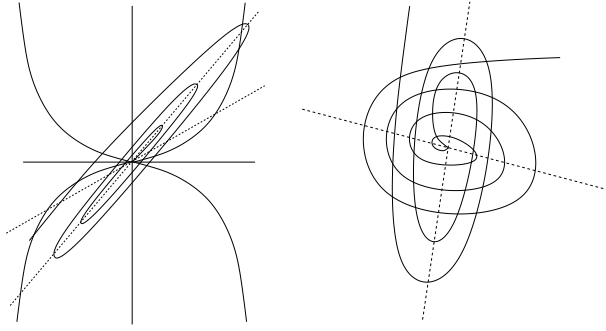}}
\caption{Tangencies with one complex foliation\label{f.onecompl}}
\end{figure}

If the pair of foliations is not of type \itt, there are two lines, intersecting in $A$, along which the foliations are tangent (see figure~\ref{f.onecompl}).

Next lemma is classical (compare~\cite[lemme 3.3]{bodu}).
\begin{lemma}
Let $(\cF,\cG)$ and $(\tilde \cF,\tilde \cG)$ be two pairs of characteristic foliations in the class \itt.

There exists a homeomorphism $h\colon U\to\tilde U$ such that $h(\cF)=\tilde \cF$ and $h(\cG)=\tilde \cG$.
\end{lemma}

So, if $\Gamma$ and $\tilde Gamma$ are in \itt, they are topologically equivalent.

Let us recall some definitions and results from~\cite{bodu}.

\subsubsection{Two logarithmic foliations}\label{ss.loglog}

Let us consider $\alpha\neq 0$, $\beta\neq 0$ ($\alpha\neq\beta$) and $\lambda > 1$ such that $t=\frac 1 2 (\lambda + \frac 1 \lambda)$ is strictly greater than $\psi(\alpha,\beta)$ (i.e. $(\cL_\alpha,\cL_\beta^\lambda)$ is not in \itt, see section~\ref{ss.bodu}). 

The foliations $\cL_\alpha$ and $\cL_\beta^\lambda$ are tangent along two lines $\Delta$ and $\Delta'$ (i.e. the vector fields $X_\alpha$ and $X_\beta^\lambda$, tangent to the foliations, are collinear). The plane is divided into two cones by $\Delta$ and $\Delta'$. 

We use the following convention to name $\Delta$ and $\Delta'$. On one cone, $S_0$, the sign of the determinant of ($X_\alpha$,$X_\beta^\lambda$) is opposite to the sign of the difference $\beta -\alpha$ (the cone $S_0$ becomes smaller as $t$ goes to $\psi(\alpha,\beta)$). We name $\Delta$ and $\Delta'$ such that the arc of the unit circle (oriented in the trigonometric way) starts on $\Delta$ and  ends on $\Delta'$.

Let us denote by $f\colon\Delta\to\Delta'$ and $g\colon\Delta\to\Delta'$ the holonomies of the restriction to $S_0$ of $\cL_\alpha$ and $\cL_\beta^\lambda$. Let $H\colon\Delta\to\Delta$ be the homothety defined by $H=g^{-1}\circ f$. And let $\mu(\alpha,\beta,t)>0$ be the ratio of this homothety. Recall that $U=\{{(\alpha,\beta,t)\in\RR^3}, {\alpha\neq 0},\ {\beta\neq 0,\ {\beta\neq\alpha}\mbox{ and }  t>\psi(\alpha,\beta)}\}$.

\begin{definition}
We set $\Psi\colon U\to\RR$ to be $\Psi(\alpha,\beta,t)=\frac \alpha {2\pi}\log(\mu(\alpha,\beta,t))$.
\end{definition}

The number $\alpha$ (resp. $\beta$) is the ratio of the homothety obtained as the holonomy of $\cL_\alpha$ (resp. $\cL_\beta^\lambda$) on $\Delta$. One can prove that $(\frac \alpha \beta, \Psi(\alpha,\beta,t))$ is a complete invariant of the pair of foliations up to topological equivalence (compare~\cite[Proposition 3.1]{bodu}). In~\cite{bodu}, we proved that $(\frac \alpha \beta, \Psi(\alpha,\beta,t))$ is a complete invariant up to the semblable relation and in particular, we proved that a pair of complex foliations with two lines of tangencies is not semblable to a pair in \itt~ (\cite[Corollaire 4.3]{bodu}). The proof works equally well for a pair of foliations with only one logarithmic foliation: the number of lines of tangencies is invariant for semblable pairs of foliations.

\subsubsection{One real and one logarithmic foliation}\label{ss.upsilon}
Let us consider $(\cL_\beta^{\lambda,\theta_0},\cR_\mu)$ a pair that is not in \itt~(see definition~\ref{d.rctt}), let $s=\frac 12(\lambda-\frac 1\lambda)$. We orient the leaves of the foliations ``from the origin". The invariant manifolds of the real foliation divide the plane into two cones, $Q_0$ and $Q_1$. As in previous section, there exist two lines of tangencies $\Delta$ and $\Delta'$. They are in the same cone $Q_0$. In this cone, each leaf of $\cR_\mu$ intersects the two lines, we denote by $\Delta$ the first and $\Delta'$ the second line intersected by an oriented leaf.

The lines $\Delta$ and $\Delta'$ divide the plane into two cones, $S_0$ and $S_1$. We define $S_0$ to be the cone included in $Q_0$. We denote by $f\colon\Delta\to\Delta'$ and $g\colon\Delta\to\Delta'$ the holonomies of the restriction to $S_0$ of $\cL_\beta^{\lambda,\theta_0}$ and $\cR_\mu$. Let $H\colon\Delta\to\Delta$ be the homothety defined by $H=g^{-1}\circ f$. And let $\nu(\beta,\mu,s,\theta_0)>0$ be the ratio of this homothety. 

\begin{definition}
We define $\Upsilon(\beta,\mu,s,\theta_0)=\frac \beta{2 \pi}\log(\nu(\beta,\mu,s,\theta_0))$ where it makes sense (on the complement of the class \itt).
\end{definition}

The number $\Upsilon(\beta,\mu,s,\theta_0)$ is a complete invariant of topological equivalence for pairs of foliations. The invariance of $\Upsilon(\beta,\mu,s,\theta_0)$ for the semblable relation is similar to the invariance of $\Psi(\alpha,\beta,t)$ in~\cite{bodu}.

\subsection{Two real foliations}\label{s.2-1.rr}
We prove here the assertions corresponding to the last row and the last column of the array of Theorem~\ref{t.2-1}. Let us consider a pair of characteristic foliations corresponding to a connection in $(2-1)^g_{\RR,\RR}$: a real foliation $\cR_p=\cR_\mu$ and $\cR_q$ the image by a diffeomorphism of a real foliation $\cR_\gamma$. We denote by $W^s(p)$, $W^{ss}(p)$, $W^u(q)$ and $W^{uu}(q)$ , the weak and strong manifolds of the origin $A$ and by $\omega^s(p)$, $\omega^{ss}(p)$, $\omega^u(q)$ and $\omega^{uu}(q)$ their tangent space at the origin.

\begin{definition}
A pair $(\cR_p,\cR_q)$ is of {\em type $(I)$} if there exists a neighbourhood $U$ of $A$ such that the foliations are transverse on $U$. Otherwise, we say that the pair is of type {\em type $(II)$}.
\end{definition}

The proof of next lemma is left to the reader.
\begin{lemma}
A pair $(\cR_p,\cR_q)$ is of {\em type $(I)$} (resp. {\em type $(II)$}) if $\omega_p^{ss}$ and $\omega_q^{uu}$ are (resp. are not) in the same connected component of $\RR P(1)\setminus\{\omega_p^s,\omega_q^u \}$.
\end{lemma}

\begin{figure}[htb]
\centerline{\includegraphics{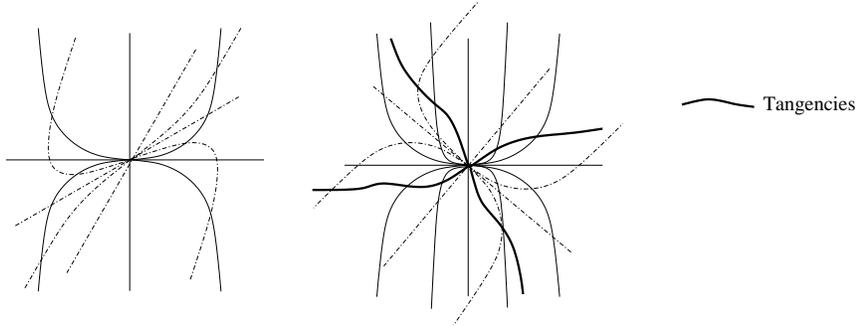}}
\caption{Type $(I)$ and $(II)$\label{f.couplesreels}}
\end{figure}

\begin{remark}\label{r.rr}
For a pair of real foliations, one can prove that two leaves of different foliations meet at most twice in a neighbourhood of the origin.

Moreover, one can prove that for a pair of type $(II)$, in a neighbourhood of the origin, there exist four segments of tangencies between the two foliations (see figure~\ref{f.couplesreels}).
\end{remark}

\begin{lemma}
If $(\cR_p,\cR_q)$ and $(\cR_{\tilde p},\cR_{\tilde q})$ have same type, there exist $U$, $\tilde U$ and $h\colon U\to\tilde U$ such that $h(\cR_p)=\cR_{\tilde p}$ and $h(\cR_q)=\cR_{\tilde q}$.
\end{lemma}
\begin{proof}
Let us deal with type $(I)$ first, we restrict the study to a neighbourhood of the origin. The invariant manifolds of $\cR_p$ divide the plane into two (smooth, not linear) cones, $C_0$ and $C_1$. One cone, namely $C_0$, contains the invariant manifolds of $\cR_q$ (therefore, it contains one cone of $\cR_q$). 

Fix $F_0$ and $G_0$ two leaves of $\cR_p$ such that $F_0$ is included in one component of $C_0\setminus\{A\}$ and $G_0$ is included in the other one. We denote by  $R_p^1$ and $R_p^2$ the parts of $F$ and $G$ which are included in the cone of $\cR_q$ included in $C_0$. Let us choose $F_1$ and $G_1$, leaves of $\cR_q$, such that $F_1$ meets one component of $C_1\setminus\{A\}$ and $G_1$ meets the other component. We denote by $R_q^1$ and $R_q^2$ the intersection of $F_1$ and $G_1$ with $C_1$. Each leaf of $\cR_p$ and $\cR_q$ meets the union $W^s(p)\cup W^u(q)\cup R_p^1\cup R_p^2 \cup R_q^1 \cup R_q^2$.

Fix any homeomorphism $h\colon W^s(p)\cup W^u(q)\cup R_p^1\cup R_p^2 \cup R_q^1 \cup R_q^2\to W^s(\tilde p)\cup W^u(\tilde q)\cup R_{\tilde p}^1\cup R_{\tilde p}^2 \cup R_{\tilde q}^1 \cup R_{\tilde q}^2$ such that $h(W^s(p))=W^s(\tilde p)$, $h(W^u(q))= W^u(\tilde q)$, \dots The homeomorphism $h$ can be extended uniquely on $U$ satisfying the conditions that $h(\cR_p)=\cR_{\tilde p}$ and $h(\cR_q)=\cR_{\tilde q}$.

For the case of type $(II)$, we remark that in a neighbourhood $U$ of the origin, a pair of foliations of type $(II)$ exhibit four curves of tangencies, with $A$ in their adherences. Let us denote by $\cD$ the union of those curves and by $W$ the union of the invariant manifolds of the pair $(\cR_p,\cR_q)$.

Fix $h\colon \cD\cup W\to \tilde \cD\cup\tilde W$ a homeomorphism mapping $A$ onto $\tilde A$ with $h(W^s(p))=W^s(\tilde p)$, $h(W^u(q))= W^u(\tilde q)$, \dots Again, the conditions $h(\cR_p)=\cR_{\tilde p}$ and $h(\cR_q)=\cR_{\tilde q}$ implie that $h$ extended uniquely to $U$.
\end{proof}

For a pair of characteristic foliations with at least one complex foliation, it is easy to find two leaves that intersect infinitely many times. So the proof of lemma~\ref{l.infini} gives the following result.
\begin{lemma}
A pair of real foliations is never semblable with a pair of characteristic foliations with at least on complex foliation.
\end{lemma} 

Following~\cite{bodu}, we define a {\em bigone} of a pair of foliations to be a disk of the plane with the origin removed, whose boundary is the union of two leaves. A {\em pearl} is a bigone such that there exists a unique curve connecting the two components of the boundary of the bigone such that the foliations are tangent along this curve. In the case of a pair of real foliations, every bigone is a pearl (compare~\cite[lemme~4.2]{bodu}).

We define a $\varepsilon-$pearl to be a ``fat" pearl: the annulus between two nested pearls. We proved in~\cite[lemme~4.8]{bodu} that if two pairs of foliations are semblable and $\Delta_\varepsilon$ is an $\varepsilon-$pearl for one pair of foliations, close enough to the origin, there exists a bigone for the second pair of foliations whose image is included in the $\varepsilon-$pearl. In particular, if two pairs of foliations are semblable, they both have (or have no) pearls.
Only type $(II)$ pairs of real foliations have pearls, so we have the next result which concludes the proof of Theorem~\ref{t.2-1}.

\begin{lemma}
If $(\cR_p,\cR_q)$ is of type $(I)$ and $(\cR_{\tilde p},\cR_{\tilde q})$ is of type $(II)$, they are not semblable.
\end{lemma} 
\bibliographystyle{amsplain}
\nocite{cela,duchaos,to} 
\bibliography{../topinv}

\end{document}